\documentclass[amsmath,english,a4paper,graphicx,12pt]{article}
\usepackage[utf8]{inputenc}
\usepackage{doc}
\usepackage{graphicx}
\usepackage{amsmath}
\usepackage{amsfonts}
\usepackage{amssymb}
\usepackage{babel}
\usepackage{latexsym}
\usepackage{mathrsfs}
\usepackage{array}

\DeclareMathOperator{\intr}{int}
\DeclareMathOperator{\sinc}{sinc}

\DeclareMathOperator{\lin}{lin}
\DeclareMathOperator{\res}{res}
\newtheorem{theorem}{Theorem}
\newtheorem{lemma}{Lemma}
\newtheorem{corollary}{Corollary}
\newtheorem{proposition}{Proposition}
\newtheorem{Remark}{Remark}
\newtheorem{definition}{Definition}

\newcommand{\C}{\mathbb{C}}
\newcommand{\R}{\mathbb{R}}
\newcommand{\N}{\mathbb{N}}
\newcommand{\li}{\left}
\newcommand{\ri}{\right}
\newcommand{\abs}[1]{\left|#1\right|}
\newcommand{\zit}[1]{\textup{(\ref{#1})}}
\newcommand{\rez}[1]{\frac{1}{{#1}}}
\newcommand{\bgl}[1]{\begin{equation}\label{#1}}
\newcommand{\egl}{\end{equation}}

\newcommand{\hal}{{\textstyle \frac{1}{2}}}

\newcommand{\Z}{\mathbb{Z}}
\newcommand{\Zz}{\mathbb{Z}\setminus \{0\}}
\newcommand{\eps}{\varepsilon}

\newcommand*{\beweis}{\begin{trivlist}\item{\it Proof.}\, }
\newcommand*{\beweisende}{\hfill$\Box$ \end{trivlist}}

\def\hs{\mathbb{H}}
\newcommand{\wT}{\widetilde{\Theta}}

\title{\bf A residue theorem for polar analytic functions and Mellin analogues of Boas' differentiation formula and Valiron's sampling formula}
\author{Carlo Bardaro, \thanks{
Department of Mathematics and Computer Sciences, University of Perugia,
via Vanvitelli 1, I-06123 Perugia, Italy, e-mail: 
carlo.bardaro@unipg.it} \and
Paul L. Butzer, \thanks{Lehrstuhl A fuer Mathematik, RWTH Aachen, Templergraben 55, Aachen, D-52056, Germany, e-mail: butzer@rwth-aachen.de}\and
Ilaria Mantellini, \thanks{Department of Mathematics and Computer Sciences, University of Perugia,
via Vanvitelli 1, I-06123 Perugia, Italy, e-mail: 
mantell@dmi.unipg.it}\and Gerhard Schmeisser\thanks{Department Mathematik, Friedrich-Alexander University Erlangen-N\"{u}rnberg, Cauerstr. 11, 91058 Erlangen, Germany, 
email: schmeisser@mi.uni-erlangen.de}
}

\date{\it Dedicated to the memory of our close Friend and Colleague Professor Domenico Candeloro, who survives in our hearts}

\begin{document}
\maketitle
\noindent

\begin{abstract}
{\small  In this paper, we continue the study of the polar analytic functions, a notion introduced in 
\cite{BBMS1} and successfully applied in Mellin analysis. Here we obtain another version of the Cauchy integral formula and a residue theorem for polar Mellin derivatives, employing the new notion of logarithmic pole. The identity theorem for polar analytic functions is also derived. As applications we obtain an analogue of Boas' differentiation formula for polar Mellin derivatives, and an extension of the classical Bernstein inequality to polar Mellin derivatives. Finally we give an analogue of the well-know Valiron sampling theorem for polar analytic functions and some its consequences.}


\end{abstract}
\section{Introduction}

In \cite{BBMS1}, we introduced  the notion of polar-analytic function as follows:
Let $\mathbb{H}$ be the half-plane $\{(r,\theta): r \in \mathbb{R}^+, \theta \in \mathbb{R}\}$ and let $f$ be a complex-valued function defined on a
neighborhood of $(r_0, \theta_0)\in \mathbb{H}.$ Then the polar derivative of $f$ at the point $(r_0, \theta_0)$, denoted by 
$(D_{{\rm pol}}f)(r_0, \theta_0),$ is given by the limit
$$\lim_{(r, \theta) \rightarrow (r_0, \theta_0)}\frac{f(r, \theta) - f(r_0, \theta_0)}{re^{i\theta} - r_0e^{i\theta_0}}.$$

Our definition leads naturally to  the classical Cauchy-Riemann equations when written in their  polar form, often treated
in the literature (\cite[Sec. 23, p. 68]{BC}; \cite[Sec. 4.3, pp. 78--82]{HI}; \cite{SI}). This notion appears as a simple way to describe functions 
which are analytic on a part of the Riemann surface of the logarithm without employing branches. The main applications of the concept of polar analyticity have been established in Mellin analysis and in the realm of quadrature formulae (see \cite{BBMS3}). In Mellin 
analysis, polar analyticity turns out to be very helpful for an efficient approach, independent of Fourier analysis. 

The corresponding notion of Mellin polar-derivative, is introduced as (see \cite{BBMS1})
$$\widetilde{\Theta}_cf(r,\theta):= re^{i\theta}(D_{{\rm pol}} f)(r, \theta) + cf(r, \theta).$$

A first development of the theory of polar analytic functions is given in the recent paper \cite{BBMS4}, in which extensions of the Cauchy integral formula and Taylor-type series are studied. 

Expansions in Fourier analysis and signal processing are often obtained by
methods of complex analysis with the residue theorem as a decisive tool.
Since in these applications only poles occur, the residue theorem can be
replaced by skillful use of Cauchy's integral formula and its extension to
derivatives. 

In the case of polar-analytic functions our analogue of Cauchy's integral
formula \cite{BBMS4} is useful in horizontal strips of width less than
$2\pi$ only since otherwise it may produce additional residues as undesired
artifacts; see \cite[formula (10)]{BBMS4}. Therefore the aim of this paper is to establish a
further analogue of Cauchy's integral formula which is always free of
artifacts, and a version of the residue theorem for polar aanalytic functions. As interesting applications, 
we establish an analogue of  Boas'  differentiation formula for polar Mellin derivatives and as a consequence, a Bernstein type inequality. Moreover, we obtain a Valiron type sampling theorem for polar analytic functions and illustrate some its consequences. In particular we establish a series representation of the polar Mellin derivative of a function $f$ belonging to  the Mellin--Bernstein space $\mathscr{B}^\infty_c$ in terms of samples of the function $f.$

We wish to dedicate the present article to our close friend and fine mathematician Professor Domenico Candeloro (1951--2019) who died prematurely on May 3. The wide range of interests of Professor Candeloro covers many aspects of mathematical analysis. In particular in his recent papers \cite{BDS},
 \cite{BDS1} he gave an interesting connection between stochastic processes, approximation theory and Mellin analysis. The last section contain a brief vita of him, as a memory of his work and life.

\section{Preliminaries}

Let $\mathbb{H}:= \{(r,\theta) \in \mathbb{R}^+ \times \mathbb{R}\}$ be the right half-plane and let ${\cal D}$ be a domain in $\mathbb{H}.$
\begin{definition}\label{def1}
We say that $f:{\cal D}\rightarrow \mathbb{C}$ is {\rm polar-analytic} on ${\cal D}$ if for any $(r_0, \theta_0) \in {\cal D}$ the limit
$$\lim_{(r,\theta) \rightarrow (r_0, \theta_0)}\frac{f(r, \theta) - f(r_0, \theta_0)}{re^{i\theta} - r_0e^{i\theta_0}} =: (D_{{\rm pol}}f)(r_0, \theta_0)$$
exists and is the same howsoever $(r, \theta)$ approaches $(r_0, \theta_0)$ within ${\cal D}.$
\end{definition}
\vskip0,3cm
The following facts are proved in \cite{BBMS4}. Polar analyticity of two functions is inherited by their arithmetic combinations and the familiar rules known for classical differentiation hold for polar derivatives as well.  Moreover, for a polar analytic function on  ${\cal D},$ $f\,:\, (r,\theta) \mapsto u(r,\theta)+ iv(r,\theta),$ where $u$ and $v$ are real-valued, identifying $\mathbb{C}$ with $\mathbb{R}^2$, we may interpret $f$
as a mapping from a subset of the half-plane $\mathbb{H}$ into $\mathbb{R}^2,$ and it turn out that $f$ is differentiable in the classical sense. Finally, 
it can be verified that $f = u + iv$ with $u,v: {\cal D}\rightarrow \mathbb{R}$ is polar-analytic on ${\cal D}$ if and only if $u$ and $v$ have continuous partial derivatives on ${\cal D}$ that satisfy the differential equations
\begin{eqnarray}\label{CRE}
\frac{\partial u}{\partial \theta} = - r \frac{\partial v}{\partial r}\,,\quad 
\frac{\partial v}{\partial \theta} = r \frac{\partial u}{\partial r}\,.
\end{eqnarray}
\vskip0,4cm
Note that these equations coincide with the Cauchy-Riemann equations of an analytic function $g$ defined by $g(z) := u(r,\theta) + i v(r, \theta)$ 
for $z= r e^{i\theta}.$ For the derivative $D_{{\rm pol}}$, we easily find that 
\begin{eqnarray}\label{dipol}
(D_{{\rm pol}}f)(r, \theta) = e^{-i\theta}\bigg(\frac{\partial}{\partial r}u(r, \theta) + i \frac{\partial}{\partial r}v (r, \theta) \bigg) = 
\frac{e^{-i\theta}}{r}\bigg(\frac{\partial}{\partial \theta}v (r, \theta) - i \frac{\partial}{\partial \theta}u (r, \theta) \bigg).
\end{eqnarray}
Since $f=u+iv,$ equations (\ref{CRE}) can be written in a more compact way as
\begin{eqnarray}\label{CRE2}
\frac{\partial f}{\partial \theta}  = ir \frac{\partial f}{\partial r}
\end{eqnarray}
and then formula (\ref{dipol}) takes the form
$$(D_{{\rm pol}}f)(r, \theta) = e^{-i\theta}\frac{\partial}{\partial r}f(r,\theta) = \frac{e^{-i\theta}}{ir} \frac{\partial}{\partial \theta}f(r,\theta).$$

Also note that $D_{{\rm pol}}$ is the ordinary differentiation on $\mathbb{R}^+.$ More precisely, if $\varphi (\cdot) := f(\cdot, 0),$ then $(D_{{\rm pol}}f)(r,0) =
\varphi'(r).$

When $g$ is an entire function, then $f: (r, \theta) \mapsto g(re^{i\theta})$ defines a function $f$ on $\mathbb{H}$ that is polar-analytic 
and $2\pi$-periodic with respect to $\theta$ and one has $(D_{{\rm pol}}f)(r, \theta) = g'(z)$ with  $z=re^{i\theta}.$ 
A converse statement is not true in general.
If $f$ is polar-analytic on $\mathbb{H}$ and  $2\pi$-periodic with respect to the second variable, 
there may not exist an entire function $h$ such that $f(r,\theta) = h(re^{i\theta}).$ A simple example is the function $f(r,\theta):= e^{-i\theta}/r.$ 
It would imply that $h(z) = 1/z$, which is analytic on $\mathbb{C}\setminus \{0\}$ only.
However, if $f$ is a polar-analytic function on $\mathbb{H},$ then $g: z=x+iy\mapsto f(e^x, y)$ is an entire function.

The main novelty of the definition of polar-analytic function is that, using this approach, we avoid periodicity with respect to the argument $\theta,$ 
and in this way we can avoid the use of Riemann surfaces. 

A simple example of a polar-analytic function that is not $2\pi$-periodic is the function $L(r, \theta):= \log r + i\theta,$ which is easily seen 
to satisfy the differential equations (\ref{CRE}). In this approach, we consider the logarithm as a single-valued function on $\mathbb{H},$ without the 
use of the Riemann surface of the logarithm $S_{{\rm log}}.$ Moreover we find ($z=re^{i\theta}$)
$$(D_{{\rm pol}} L)(r, \theta) = e^{-i\theta}\frac{1}{r} = \frac{1}{re^{i\theta}} = \frac{1}{z}\,.$$

As we anticipated in the introduction, for a fixed real number $c \in \mathbb{R},$ we define the polar derivative of a polar-analytic function $f$ in Mellin frame by the formula
\begin{eqnarray}\label{mellinpolar}
\widetilde{\Theta}_cf(r,\theta):= re^{i\theta}(D_{{\rm pol}} f)(r, \theta) + cf(r, \theta),
\end{eqnarray}
provided that the polar derivative $D_{{\rm pol}} f$ exists at the point $(r,\theta) \in \mathbb{H}.$ As before, for $\varphi(r):= f(r, 0),$ we have $\widetilde{\Theta}_cf(r,0) = (\Theta_c\varphi)(r),$ where $\Theta_c$ is the usual Mellin differential operator (see \cite{BJ}).

The higher order polar Mellin derivatives may be defined through the representation formula for (usual) Mellin derivatives in terms of Stirling numbers of the second type (see \cite[Lemma 9]{BJ}), namely ($k \in \mathbb{N}$),
\begin{eqnarray}\label{higherorder}
\widetilde{\Theta}^{k}_cf(r,\theta):= \sum_{j=0}^{k} S_c(k,j)r^je^{ij\theta} D_{\rm pol}^{j}f(r,\theta) \quad \quad ((r,\theta) \in \mathbb{H}).
\end{eqnarray}

\section{A general Cauchy's integral formula for polar analytic functions}

First we recall here the following results (see \cite{BBMS4}).
\begin{lemma}\label{lemma1}
Let $\mathcal{D}$ be a convex domain in $\mathbb{H}.$ Let $g:\mathcal{D}\rightarrow \mathbb{C}$ be continuous on the whole of $\mathcal{D}$ and polar-analytic except at a point $(r_0,\theta_0) \in \mathcal{D}.$ Suppose that $\gamma$ is a positively oriented, closed, regular curve in $\mathcal{D}$ that is the boundary of a convex domain $\mbox{{\rm int}}(\gamma).$ Then
$$\int_\gamma g(r,\theta)e^{i\theta}(dr + ird\theta) = 0.$$
\end{lemma}

\begin{theorem}\label{CIF}
Let $\mathcal{D}$ be a convex domain in $\mathbb{H},$ and let $f:\mathcal{D}\rightarrow \mathbb{C}$ be polar-analytic on $\mathcal{D}.$ Let 
$\gamma$ be a positively oriented, closed, regular curve that is the boundary of a convex domain $\intr(\gamma) \subset \mathcal{D}.$ 
Given $(r_0, \theta_0) \in \intr(\gamma),$ define $\theta_j:= \theta_0 + 2j\pi$ for $j \in \mathbb{Z}.$ Suppose that none of the points 
$(r_0,\theta_j)$ lies on $\gamma.$ Then
\begin{eqnarray}\label{cauchyformula}
\frac{1}{2\pi i}\int_\gamma \frac{f(r,\theta) e^{i\theta}}{re^{i\theta} - r_0 e^{i\theta_0}}(dr + ir d\theta) = 
\sum_{(r_0, \theta_j)\in \intr(\gamma)} f(r_0, \theta_j).
\end{eqnarray}
In particular, if $\gamma$ lies in a strip $\mathbb{H}_{a,b}$ with $0<b-a<2\pi,$ then
\begin{eqnarray}\label{strip}
\frac{1}{2\pi i}\int_\gamma \frac{f(r,\theta) e^{i\theta}}{re^{i\theta} - r_0 e^{i\theta_0}}(dr + ir d\theta) = f(r_0, \theta_0).
\end{eqnarray}
\end{theorem}
\vskip0,4cm
We begin with the following extension of \cite[Theorem 6.1]{BBMS4}.
\begin{theorem}\label{thm3}
Let $f\,:\,\mathcal{D} \to\C$ be polar-analytic on $\mathcal{D}$ and let
$c\in\R$. If $(r_0,\theta_0)\in\mathcal{D}$, then there holds the
expansion
\bgl{thm31}
\bigl(re^{i\theta}\bigr)^c f(r,\theta)\,=\,
\bigl(r_0e^{i\theta_0}\bigr)^c \sum_{k=0}^\infty
\frac{\bigl(\wT_c^kf\bigr)(r_0,\theta_0)}{k!} \li(\log \frac{r}{r_0} +
i(\theta-\theta_0)\ri)^k,
\egl
converging uniformly on every polar-disk
$E\bigl((r_0,\theta_0),\rho\bigr)\subset \mathcal{D}$.
\end{theorem}

\beweis
We only have to extend the proof of \cite[Theorem 6.1]{BBMS4} slightly. Let
$$A\,:=\,\li\{z=x+iy\in\C\::\: (e^x,y)\in\mathcal{D}\ri\}.$$
In \cite{BBMS4} we have shown that
$$ g\::\: z=x+iy \, \longmapsto\, f(e^x,y)$$
is analytic on $A$ and for $z_0:=\log r_0+i\theta_0$, we have
$$ g^{(k)}(z_0)\,=\,\bigl(\wT_0^kf\bigr)(r_0,\theta_0).$$
Now define $h(z):= e^{cz} g(z)$. Then $h$ is also analytic on $A$ and
\begin{align*}
h^{(k)}(z_0) &= e^{cz_0} \sum_{j=0}^k {k\choose j} g^{(j)}(z_0)
c^{k-j}\\[1ex]
&= \bigl(r_0e^{i\theta_0}\bigr)^c \sum_{j=0}^k {k\choose j}
\bigl(\wT_0^jf\bigr)(r_0,\theta_0) c^{k-j}.
\end{align*}
It can be shown by induction on $k$ that the sum on the right-hand side is
equal to $\bigl(\wT_c^kf\bigr)(r_0,\theta_0).$ Since
by Taylor expansion
$$ e^{cz} g(z)\,=\, \sum_{k=0}^\infty \frac{h^{(k)}(z_0)}{k!}\,(z-z_0)^k,$$
we obtain \zit{thm31} by substituting $z=x+iy$ with $x=\log r$ and
$y=\theta$. The statement on convergence is seen as in the \cite{BBMS4}.
\beweisende
\vskip0,3cm
As a consequence of the series expansion (\ref{thm31}) we deduce the following identity theorem for polar analytic functions.
\begin{theorem}\label{thm1i}
Let $\mathcal{D}$ be a domain in $\hs$ and let $f\,:\,\mathcal{D} \to \C$ be polar-analytic.
Suppose that $(r_0, \theta_0)\in\mathcal{D}$ is an accumulation point of distinct zeros of $f$.
Then $f$ is identically zero.
\end{theorem}

For the proof we make decisive use of an observation which we state as a lemma for easy reference.

\begin{lemma}\label{lem1i}
Under the hypotheses of Theorem~\ref{thm1i}, it follows that
$$\bigl(\wT_0^kf\bigr)(r_0, \theta_0)=0 \quad \hbox{ for all } k\in\N_0.$$
\end{lemma}

\beweis
For a proof by contradiction, we assume that there exists a smallest non-negative integer $k_0$ such that
$\bigl(\wT_0^{k_0}f\bigr)(r_0, \theta_0)\ne 0$. Then $f$ has an expansion
\begin{align*}
f(r,\theta)=&\li(\log \frac{r}{r_0} +i(\theta-\theta_0)\ri)^{k_0}\\[1ex]
&\times\li(\frac{\bigl(\wT_0^{k_0}f\bigr)(r_0,\theta_0)}{k_0!}
+ \sum_{k=k_0+1}^\infty\frac{\bigl(\wT_0^kf\bigr)(r_0,\theta_0)}{k!}\li(\log \frac{r}{r_0}+i(\theta-\theta_0)\ri)^{k-k_0}\ri)
\end{align*}
converging in a polar disk $E((r_0,\theta_0),\rho)$ with $\rho>0$ by Theorem \ref{thm3}; see also \cite[Theorem~6.1]{BBMS4}.
On the right-hand side, the first factor vanishes for $(r,\theta)=(r_0,\theta_0)$ only if $k_0\ge 1$; otherwise it is identically $1$.
The second factor is different from zero at $(r,\theta)=(r_0,\theta_0)$. By continuity, it is also different
from zero in a sufficiently small neighborhood of $(r_0, \theta_0)$. Hence $(r_0, \theta_0)$ cannot be an
accumulation point of distinct zeros of $f$. A contradiction!
\beweisende

\begin{trivlist}\item{\it Proof of Theorem~\ref{thm1i}.}\,
Let $(r^\ast, \theta^\ast)$ be an arbitrary point in $\mathcal{D}$. It suffices to show that $f(r^\ast, \theta^\ast)=0.$
Let $\gamma\,:\,[0,1] \to \mathcal{D}$ be a Jordan arc in $\mathcal{D}$ such that $\gamma(0)=(r_0, \theta_0)$ and $\gamma(1)=(r^\ast, \theta^\ast)$.
Since the trace of $\gamma$ has a positive distance from the boundary of $\mathcal{D}$, there exists a number $\rho_{\rm inf}>0$
such that for each point $(r, \theta)$ on $\gamma$ there is a polar disk $E((r, \theta),\rho)\subset \mathcal{D}$ with $\rho\ge\rho_{\rm inf}$.

First consider  a polar disk $E_0:=E((r_0, \theta_0), \rho_0)\subset\mathcal{D}$ with  $\rho_0\ge\rho_{\rm inf}.$
By Theorem \ref{thm3} we have a representation
$$f(r, \theta)\,=\, \sum_{k=0}^\infty \frac{\bigl(\wT_0^kf\bigr)(r_0, \theta_0)}{k!}\li(\log \frac{r}{r_0} +i(\theta-\theta_0)\ri)^k$$
holding for all $(r,\theta)\in E_0$.
Now, by Lemma~\ref{lem1i}, it follows that the restriction of $f$ to $E_0$ is identically zero. Thus, if $(r^\ast, \theta^\ast)\in E_0$, we have
reached our aim.

Otherwise, let $(r_1, \theta_1):=\gamma(t_1)$ with $0<t_1\le1$ be a point where $\gamma$ intersects the boundary of $E_0$.
By our previous observations, there exists a polar disk $E_1:=E((r_1, \theta_1), \rho_1)\subset \mathcal{D}$ with $\rho_1\ge\rho_{\rm inf}$
such that
$$f(r, \theta)\,=\, \sum_{k=0}^\infty \frac{\bigl(\wT_0^kf\bigr)(r_1, \theta_1)}{k!}\li(\log \frac{r}{r_1} +i(\theta-\theta_1)\ri)^k$$
for all $(r,\theta)\in E_1$. Since $(r_1, \theta_1)$, being a boundary point of $E_0$, is an accumulation point of distinct zeros of $f$,
we conclude with the help of Lemma~\ref{lem1i} that the restriction of $f$ to $E_1$ vanishes identically. If $(r^\ast, \theta^\ast)\not\in E_1$,
we continue this procedure with a point $(r_2, \theta_2):=\gamma(t_2)$, where $t_1< t_2 \le1$, lying on the boundary of $E_1$.
Since $\rho_{\rm inf}>0$, we arrive after a finite number of steps at a polar disk $E_n$, say, such that the restriction of $f$ to $E_n$
vanishes identically and $(r^\ast, \theta^\ast)\in E_n$. This completes the proof.
\beweisende

\vskip0,4cm
Now, we state the announced Cauchy's integral formula for polar analiytic functions in its general form.

\vskip0,2cm
\begin{theorem}\label{thm4}
Let $\mathcal{D}$ be a convex domain in $\hs$, and let $f\,:\, \mathcal{D}
\to \C$ be polar-analytic on $\mathcal{D}$. Let $\gamma$ be a positively
oriented, closed, regular curve that is the boundary of a convex domain
$\intr(\gamma)\subset \mathcal{D}$. Then, for $(r_0, \theta_0)
\in\intr(\gamma)$, $c\in\R$ and $k\in\N_0$, we have
\bgl{thm41}
\rez{2\pi i} \int_\gamma \frac{(re^{i\theta})^{c-1} f(r,\theta)
e^{i\theta}}{\bigl(\log(r/r_0)+i (\theta-\theta_0)\bigr)^{k+1}} 
(dr +ir d\theta)\,=\,
\bigl(r_0e^{i\theta_0}\bigr)^c 
\frac{(\widetilde{\Theta}^k_cf)(r_0, \theta_0)}{k!}\,.
\egl
\end{theorem}

\beweis
We decompose $\intr(\gamma)$ into four parts as in Figure \ref{fig1} (see also \cite{BBMS4}). 
Further we recall that the disk is of radius
$\eps>0$ and has a circle $\lambda_\eps$ as boundary, which shall be 
positively oriented this time.
\begin{figure}[htbp]
\centering
\includegraphics[width=70mm]{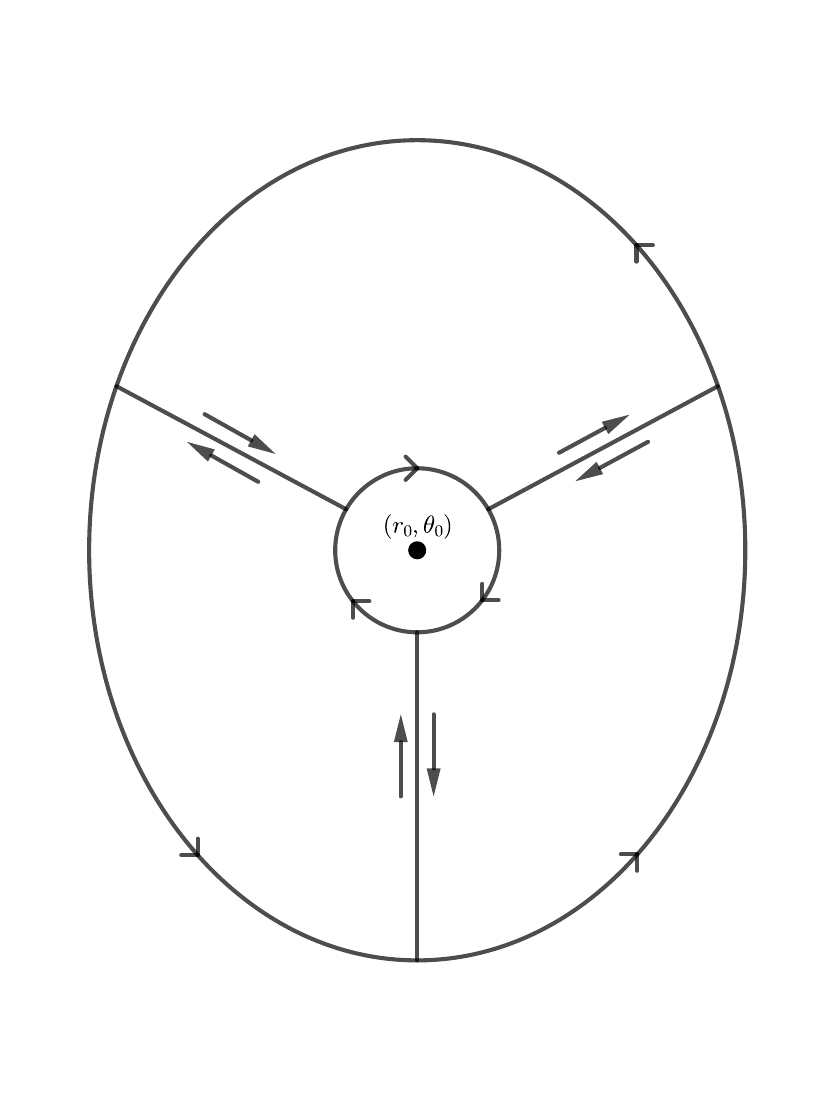}%
\caption{\small Decomposition of ${\rm int}(\gamma)$ into four parts}\label{fig1} 
\end{figure}

Furthermore $\beta_1$, $\beta_2$ and $\beta_3$ denote the positively oriented
boundaries of the other three parts. Note that each $\beta_j$ lies in a
convex subset of $\mathcal{D}$ where the function
\bgl{thm4p1}
F\::\: (r, \theta)\,\longmapsto\, \frac{(re^{i\theta})^{c-1}
f(r,\theta)}{\bigl(\log(r/r_0)+i(\theta-\theta_0)\bigr)^{k+1}}
\egl
is polar-analytic. By \cite[Theorem 4.1]{BBMS4} we have
$$ \int_{\beta_j} F(r,\theta) e^{i\theta} (dr+ ir d\theta)\,=\,0
\quad (j=1,2,3).$$
This enables us to conclude that
\bgl{thm4p2}
\rez{2\pi i} \int_\gamma F(r,\theta) e^{i\theta} (dr + ir d\theta)\,=\,
\rez{2\pi i} \int_{\lambda_\eps} F(r,\theta) e^{i\theta}(dr + ir
d\theta).
\egl
From Theorem~\ref{thm3} we know that
$$F(r,\theta)\,=\,\bigl(r_0e^{i\theta_0}\bigr)^{c-1} \sum_{\ell=0}^\infty
\frac{\bigl(\wT_{c-1}^\ell f\bigr)(r_0,\theta_0)}{\ell!}
\li(\log \frac{r}{r_0} + i(\theta-\theta_0)\ri)^{\ell-k-1}.$$
Note that on the right-hand side all terms with $\ell\ge k+1$ are
polar-analytic as functions of $(r, \theta)$. Thus
\bgl{thm4p3}
\rez{2\pi i} \int_{\lambda_\eps} F(r,\theta) e^{i\theta} (dr+ ir d\theta)
\,=\, \bigl(r_0 e^{i\theta_0}\bigr)^{c-1} \sum_{\ell=0}^k
\frac{\bigl(\wT_{c-1}^\ell f\bigr)(r_0, \theta_0)}{\ell!} I_{k+1-\ell},
\egl
where
$$ I_j\,:=\, \rez{2\pi i} \int_{\lambda_\eps}\frac{e^{i\theta} (dr+
ird\theta)}{\bigr(\log(r/r_0) + i(\theta-\theta_0)\bigr)^j}
\qquad (j=1, \dots, k+1).$$
Now consider the mapping $(r,\theta) \mapsto z=re^{i\theta}$. If $(r,
\theta)$ traverses $\lambda_\eps$ once, then, for sufficiently small
$\eps>0$, the point $z$ runs along a piecewise continuously
differentiable, closed curve $\gamma_\eps$, say, in $\C$, surrounding
$z_0=r_0e^{i\theta_0}$ exactly once in the mathematically positive sense.
Therefore,
$$I_j = \rez{2\pi i} \int_{\gamma_\eps} \frac{dz}{(\log z - \log z_0)^j}
= \rez{2\pi i} \int_{\log\circ\gamma_\eps}\frac{e^w dw}{(w-\log z_0)^j}
=\frac{r_0e^{i\theta_0}}{(j-1)!}\,.$$
With this result \zit{thm4p3} becomes
\begin{align*}
\rez{2\pi i} \int_{\lambda_\eps} F(r,\theta) e^{i\theta} (dr+ ir d\theta)
&= \bigl(r_0 e^{i\theta_0}\bigr)^c \sum_{\ell=0}^k
\frac{\bigl(\wT_{c-1}^\ell f\bigr)(r_0, \theta_0)}{\ell! (k-\ell)!} \\[1ex]
&=\frac{\bigl(r_0 e^{i\theta_0}\bigr)^c}{k!} \sum_{\ell=0}^k 
{k\choose \ell} \bigl(\wT_{c-1}^\ell f\bigr)(r_0,\theta_0)\\[1ex]
&=  \frac{\bigl(r_0 e^{i\theta_0}\bigr)^c}{k!} \bigl(\wT_c^kf\bigr)(r_0,
\theta_0).
\end{align*}
The last conclusion is seen by noting that $\wT_cf=\wT_{c-1}f +f$ and
using this relation repeatedly. This completes the proof.
\beweisende
\bigskip

\section{A residue theorem for polar-analytic functions}

First we introduce some notions.

\begin{definition}\label{def1}
Let $(r_0,\theta_0)\in\hs$ and let $\mathcal{U}\subset\hs$ be an open
neighborhood of $(r_0,\theta_0)$.
\begin{enumerate}
\item[{\rm(i)}] If $f\,:\,\mathcal{U}\setminus\bigl\{(r_0,\theta_0)\bigr\}
\to\C$ is polar-analytic, then $(r_0,\theta_0)$ will be called an
{\em isolated singularity} of $f$.
\item[{\rm(ii)}]
An isolated singularity $(r_0,\theta_0)$ is said to be a {\em logarithmic
pole of order} $k$ if $k\in\N$ and there exists a polar-analytic function 
$g\,:\,\mathcal{U}\to\C$ with $g(r_0,\theta_0)\ne0$ such that
$$ f(r,\theta)\,=\,
\frac{g(r,\theta)}{\bigl(\log(r/r_0)+i(\theta-\theta_0)\bigr)^k}
\quad \hbox{ for }
(r,\theta)\in\mathcal{U}\setminus\bigl\{(r_0,\theta_0)\bigr\}.$$
In this case
$$\bigl(\res_cf\bigr)(r_0,\theta_0)\,:=\, \bigl(r_0e^{i\theta_0}\bigr)^c
\, \frac{\bigl(\wT_c^{k-1}g\bigr)(r_0, \theta_0)}{(k-1)!}$$
will be called the $c$-{\em residue} of $f$ at $(r_0,\theta_0)$.
\end{enumerate}
\end{definition}

We are now in a position to formulate a residue theorem for logarithmic
poles which will be suitable for many applications in Mellin analysis.

\begin{theorem}\label{thm5}
Let $\mathcal{D}$ be a convex domain in $\hs$ and let $f$ be
polar-analytic on $\mathcal{D}$ except for isolated singularities which
are all logarithmic poles. Let $\gamma$ be a positively oriented, closed,
regular curve that is the boundary of a convex domain
$\intr(\gamma)\subset\mathcal{D}$. Suppose that no isolated singularity
lies on $\gamma$ while $(r_j, \theta_j)$ for $j=1, \dots, m$ are the
isolated singularities lying in $\intr(\gamma).$ Then, for $c\in\R$, there
holds
$$
\int_\gamma \bigl(re^{i\theta}\bigr)^{c-1}
f(r,\theta)e^{i\theta}\,(dr+ird\theta)\,=\, 2\pi i \sum_{j=1}^m
\bigl(\res_cf\bigr)(r_j,\theta_j).
$$
\end{theorem}

\beweis
For $(x_0, y_0)\in\R^2$ and $\eps>0$, we consider the vertical grid
$$V(x_0, \eps)\,:=\, \bigl\{(r,\theta)\in\R^2\::\: \theta\in\R,\, 
r=x_0+n\eps,\, n\in\Z\bigr\}$$
and the horizontal grid
$$H(y_0, \eps)\,:=\, \bigl\{(r,\theta)\in\R^2\::\: r\in\R, \,
\theta=y_0+n\eps,\, n\in\Z\bigr\}.$$
Their union constitutes a net
$$ N\bigl((x_0,y_0), \eps\bigr)\,:=\, V(x_0, \eps) \cup H(y_0,\eps).$$
The intersection $\intr(\gamma)\cap N\bigl((x_0,y_0), \eps\bigr)$
creates a tessellation of $\intr(\gamma)$ into squares with edges of
length $\eps$ and further convex sets whose boundary contains a piece of
$\gamma$. The latter may be called boundary sets. The collection of all
subsets of the tessellation shall be denoted by $\mathcal{T}$. If $\eps>0$
is sufficiently small, then the boundary sets will be free of isolated
singularities and different isolated singularites in $\intr(\gamma)$ 
will lie in different squares of the tessellation. Furthermore, if
$(r_j, \theta_j)$ lies in the interior of a square, then it will remains
there under all sufficiently small variations of $(x_0,y_0)$. But if
$(r_j, \theta_j)$ lies on the boundary of a square, then there exist
arbitrarily small variations such that $(r_j, \theta_j)$ goes inside.
We may therefore assume that each $(r_j, \theta_j)$ lies inside a square
$Q_j \in\mathcal{T}$ for $j=1, \dots, m$. For each set
$P\in\mathcal{T}$, we denote by $\partial P$ its
positively oriented boundary. Note that each line segment in
$\intr(\gamma)$ which comes from a mesh of the net belongs to the 
boundaries of exactly
two subsets of the tessellation where it occurs with opposite orientations.
Therefore
\begin{align*}
\int_\gamma \bigl(re^{i\theta}\bigr)^{c-1} f(r,\theta) e^{i\theta}
(dr +ird\theta) &=
\sum_{P\in\mathcal{T}} \int_{\partial P} \bigl(re^{i\theta}\bigr)^{c-1}
f(r, \theta) e^{i\theta} (dr+ir d\theta)\\[1ex]
&= \sum_{j=1}^m \int_{\partial Q_j}\bigl(re^{i\theta}\bigr)^{c-1} 
f(r,\theta) e^{i\theta}(dr + ir d\theta)
\end{align*}
since by \cite[Theorem~4.1]{BBMS4} the integral along $\partial P$ vanishes if
$P$ does not contain an isolated singularity. 
As all our isolated singularities are logarithmic poles, we conclude with
the help of Theorem~\ref{thm4} and by using the notion of residue that
$$
\int_{\partial Q_j}\bigl(re^{i\theta}\bigr)^{c-1} 
f(r,\theta) e^{i\theta}(dr + ir d\theta)\,=\, 2\pi i
\bigl(\res_cf\bigr)(r_j,\theta_j) \quad (j=1, \dots, m).$$
This completes the proof.
\beweisende
\bigskip

\section{Boas'  differentiation formula
for polar Mellin derivatives}

We recall that (see \cite{BBMS1}) the Mellin--Bernstein space $\mathscr{B}^p_{c,T}$ comprises
all functions $f\,:\, \hs \to\C$ with the following properties:
\begin{enumerate}
\item[(i)] \, $f$ is polar-analytic on $\hs$;
\item[(ii)] \, $f(\cdot ,0) \in X_c^p$;
\item[(iii)] \, there exists a constant $C_f>0$ such that
$r^c \abs{f(r, \theta)} \le C_f e^{T\abs{\theta}}$ for all 
$(r, \theta)\in \hs.$
\end{enumerate}

We now state three useful assertions on transformations in Mellin--Bernstein
spaces. They are verified by straightforward calculations. To show polar
analyticity, we simply check that the Cauchy--Riemann
equations in polar form are satisfied. The first statement affirms that
$\mathscr{B}^p_{c,T}$ is invariant under Mellin translations.

\begin{proposition}\label{prop1}
Let $f\in \mathscr{B}^p_{c,T}$, where $p\ge 1$, $c\in\R$ and $T>0$. 
For $t>0$, define
$$ g\::\: (r, \theta) \, \longmapsto \, t^c\,f(tr, \theta).$$
Then $g\in\mathscr{B}^p_{c,T};$ in particular, $\|g(\cdot, 0)\|_{X_c^p}
=\|f(\cdot, 0)\|_{X_c^p}$ and
$r^c \abs{g(r,\theta)}\le C_f e^{T\abs{\theta}}$ for all
$(r,\theta)\in\hs$. Furthermore,
$$\bigl(\widetilde{\Theta}_cg\bigr)(r,\theta)\,=\,
t^c\bigl(\widetilde{\Theta}_cf\bigr)(tr,\theta).$$
\end{proposition}

The second statement concerns a transformation of $T$.

\begin{proposition}\label{prop2}
Let $f\in \mathscr{B}^p_{c,T}$, where $p\ge 1$, $c\in\R$ and $T>0$. Define
$$ h\::\: (r, \theta) \, \longmapsto \, f(r^{1/T}, \theta/T).$$
Then $h\in\mathscr{B}^p_{c/T,1};$ in particular, $\|h(\cdot, 0)\|_{X_{c/T}^p}
=\|f(\cdot, 0)\|_{X_c^p}$ and
$r^{c/T} \abs{h(r,\theta)}\le C_f e^{\abs{\theta}}$ for all
$(r,\theta)\in\hs$.
Furthermore,
$$\bigl(\widetilde{\Theta}_{c/T}h\bigr)(r,\theta)\,=\, \rez{T}\,
\bigl(\widetilde{\Theta}_cf\bigr)\li(r^{1/T}, \frac{\theta}{T}\ri).$$
\end{proposition}

The third statement concerns a shift of the second argument of $f$.
For verifying property (ii) in the definition of Mellin--Bernstein spaces  we will use Theorem 4.1 in \cite{BBMS1}.
\begin{proposition}\label{prop3}
Let $f\in\mathscr{B}^p_{c,T}$, where $p\ge 1$, $c\in\R$ and $T>0$.
For $\alpha\in\R$, define
$$ \phi\::\: (r,\theta)\,\longmapsto\, f(r, \theta+\alpha).$$
Then $\phi\in\mathscr{B}^p_{c,T}$; in particular,
$\|\phi(\cdot,0)\|_{X_c^p}\le e^{T\abs{\alpha}}\,\|f(\cdot,0)\|_{X_c^p}$
and $r^c \abs{\phi(r, \theta)}\le C_f e^{T(\abs{\alpha}+\abs{\theta})}$
for all $(r, \theta)\in \hs$. Furthermore,
$$ \bigl(\widetilde{\Theta}_c\phi\bigr)(r,\theta)\,=\, 
\bigl(\widetilde{\Theta}_cf\bigr)(r, \theta+\alpha).$$
\end{proposition}

Now we present a formula for the polar Mellin derivatives. It is analogous
to a differentiation formula of Boas for bandlimited functions.

\begin{theorem}\label{thm2}
Let $f\in \mathscr{B}^p_{c,T}$, where $p\ge 1$, $c\in\R$ and $T>0$.
Then
\bgl{thm21}
\bigl(\widetilde{\Theta}_cf\bigr)(r, \theta)\,=\,
\frac{4T}{\pi^2} \sum_{k\in\Z} \frac{(-1)^k}{(2k+1)^2}\,\,
e^{(k+1/2)\pi c/T} f\bigl(r e^{(k+1/2)\pi/T}, \theta\bigr)
\egl
for $(r,\theta)\in\hs$. Multiplied by $r^c$,
the series converges absolutely and uniformly on strips of bounded width
parallel to the $r$-axis in $\hs$.
\end{theorem}

\beweis
For simplicity, we first suppose that $T=1$. Consider the function
$$ F\::\, (r,\theta) \, \longmapsto \,
\frac{f(r,\theta)}%
{(\log r+i\theta)^2 \cos(\log r+i\theta)}\,.$$
It is polar-analytic on $\hs$ except for isolated singularities at the points where the denominator
vanishes. Writing $r_k:= e^{(k+1/2)\pi}$ for short, we obtain the
exceptional points as $(1,0)$ and $(r_k, 0)$ for $k\in\Z$.
The first one is a logarithmic pole of order two; all the others
are logarithmic poles of order one.

Next, for $n\in\N$, let $\mathcal{R}_n$ be the rectangle with vertices at
$(e^{\pm n\pi}, \pm n\pi)$ and denote by $\partial\mathcal{R}_n$ its
positively oriented boundary. We decompose $\mathcal{R}_n$ into $2n+1$
rectangular parts by vertical lines passing through the points
$(e^{\pm \pi/4}, 0)$ and $(e^{k\pi}, 0)$ for $k= \pm1, \pm2, \dots,
\pm(n-1).$ Note that each part contains exactly one exceptional point
in its interior. Denote by $\mathcal{P}$ the part that contains $(1, 0)$
and by $\mathcal{Q}_k$ the one that contains $(r_k, 0)$ for
$k= -n, -n+1, \dots , n-1$. Let $\partial\mathcal{P}$ and
$\partial\mathcal{Q}_k$ be the positively oriented boundaries of
$\mathcal{P}$ and $\mathcal{Q}_k$, respectively.
\begin{figure}[htbp]
\includegraphics[width=180mm]{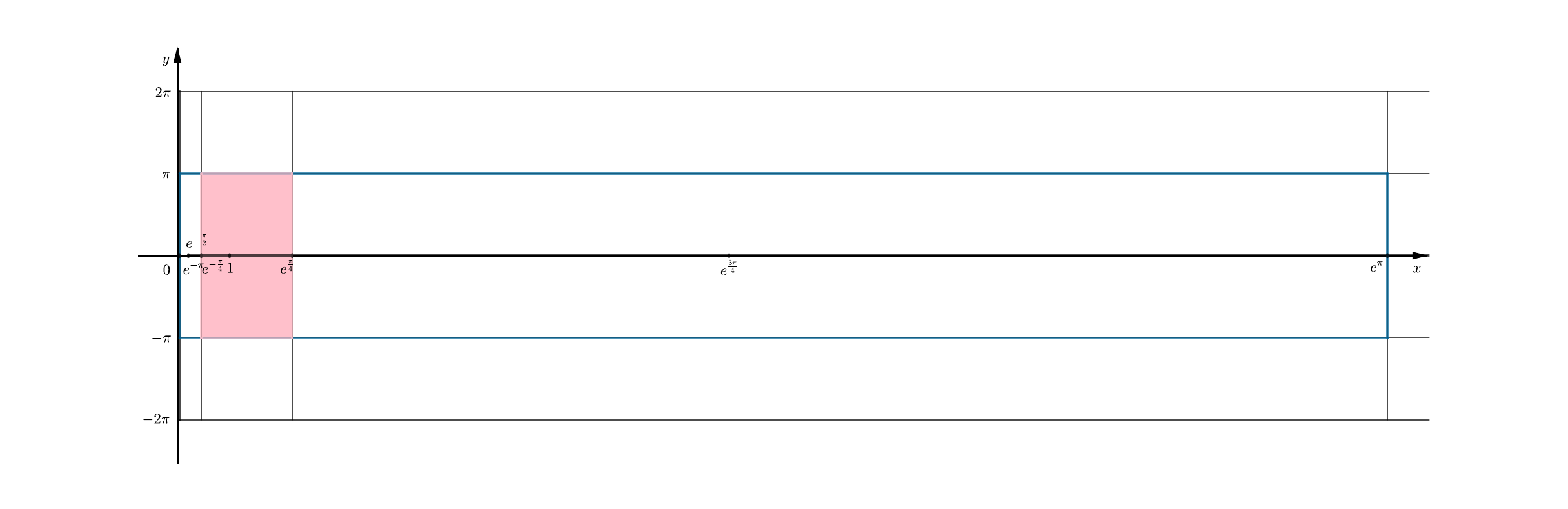}%
\caption{\small Decomposition of the rectangle $\mathcal{R}_1.$ The region $\mathcal{P}$ is colored.}\label{fig2} 
\end{figure}
Employing Theorem~\ref{thm5}, we find that
\bgl{thm2p1a}
\rez{2\pi i} \int_{\partial\mathcal{R}_n}
\bigl(re^{i\theta}\bigr)^{c-1}F(r,\theta)e^{i\theta} (dr +ird\theta)
= \bigl(\res_cF\bigr)(1,0)
 + \sum_{k=-n}^{n-1} \bigl(\res_cF\bigr)(r_k,0).
\egl
Introducing
$$ \phi(r,\theta)\,:=\, \frac{f(r, \theta)}{\cos(\log r +i\theta)}\,,$$
we have
$$ \bigl(\res_cF\bigr)(1,0)\,=\,\bigl(\wT_c\phi\bigr)(1,0)\,=\,
\bigl(\wT_cf\bigr)(1,0).$$
For calculating the other $c$-residues, we factor the cosine as
$$\cos(\log r + i\theta)\,=\, \bigl(\log r+i\theta - (k+\hal)\pi\bigr) 
\,\psi_k(r,\theta),$$
where $\psi_k$ is polar-analytic in a neighborhood of $(r_k, 0)$ and
$$ \psi_k(r_k,0)\,=\, \lim_{(r,\theta)\to (r_k,0)}
\frac{\cos(\log r+i\theta)}{\log r+i\theta-(k+\hal)
\pi} \,=\, (-1)^{k+1}.$$
Thus,
$$ \bigl(\res_cF\bigr)(r_k,0)\,=\,
\frac{r_k^c f(r_k,0)}{\bigl((k+\hal)\pi\bigr)^2 \psi(r_k,0)}
\,=\, \frac{4}{\pi^2}\,\frac{(-1)^{k+1}}{(2k+1)^2} \,r_k^c f(r_k,0).
$$
With these values of the $c$-residues, we may rewrite
\zit{thm2p1a} as
\begin{align*}
\bigl(\widetilde{\Theta}_c f\bigr)(1,0) &=
\frac{4}{\pi^2} \sum_{k=-n}^{n-1} \frac{(-1)^k}{(2k+1)^2} r_k^c
f(r_k,0)\\[1ex]
&\quad + \rez{2\pi i} \int_{\partial\mathcal{R}_n} 
\bigl(re^{i\theta}\bigr)^{c-1}F(r,\theta) e^{i\theta}
(dr + ir d\theta).
\end{align*}

Next we want to estimate the integral on the right-hand side. Noting that
$$\abs{\cos(\log r +i\theta)} > \frac{e^{\abs{\theta}}}{3}$$
on $\partial\mathcal{R}_n$ and employing properity (iii) of functions
belonging to $\mathscr{B}^p_{c,1}$, we see that
$$ r^{c-1}\abs{F(r,\theta)}\,\le\, \frac{3C_f}{r(\log^2r +\theta^2)}
\quad \hbox{ for } (r,\theta)\in \partial\mathcal{R}_n.$$
Considering the integrals along the vertical and the horizontal line
segments of $\partial\mathcal{R}_n$ separately, we find that
\begin{align*}
& \abs{\rez{2\pi i} \int_{\partial\mathcal{R}_n} \bigl(re^{i\theta}\bigr)^{c-1}
F(r,\theta) e^{i\theta} (dr +ir d\theta)}\\[1ex]
&\qquad\qquad \le \frac{3C_f}{\pi}\li( \int_{-n\pi}^{n\pi}\frac{d\theta}%
{(n\pi)^2 + \theta^2} + \int_{e^{-n\pi}}^{e^{n\pi}}\rez{\log^2 r
+(n\pi)^2}\,\frac{dr}{r}\ri)\\[1ex]
&\qquad\qquad =\frac{6C_f}{n\pi^2} \int_{-1}^1 \frac{dx}{1+x^2}
\:\longrightarrow\, 0 \hbox{ as } n\to\infty.
\end{align*}
Thus
$$
\bigl(\widetilde{\Theta}_c f\bigr)(1,0) \,=\,
\frac{4}{\pi^2} \sum_{k=-\infty}^\infty \frac{(-1)^k}{(2k+1)^2}\, r_k^c
f(r_k,0)\quad \hbox{ for } f\in\mathscr{B}^p_{c,1}.$$

With the Mellin translation of Proposition~\ref{prop1}, we deduce
$$
\bigl(\widetilde{\Theta}_c f\bigr)(t,0) \,=\,
\frac{4}{\pi^2} \sum_{k=-\infty}^\infty \frac{(-1)^k}{(2k+1)^2}\, r_k^c
f(tr_k,0)\quad (t>0),$$
valid for $f\in\mathscr{B}^p_{c,1}.$

Next we want to extend this formula to $f \in\mathscr{B}^p_{c,T}$
with an arbitrary $T>0$. For this, we note that the formula is valid for
the function $h$ defined in Proposition~\ref{prop2} provided that we
replace $c$ with $c/T$. Expressing $h$ in terms of $f$, we obtain
$$
\bigl(\widetilde{\Theta}_c f\bigr)(t^{1/T},0) \,=\,
\frac{4T}{\pi^2} \sum_{k=-\infty}^\infty \frac{(-1)^k}{(2k+1)^2}\,\, 
r_k^{c/T} f(t^{1/T}r_k^{1/T},0)\quad (t>0),$$
now valid for $f \in\mathscr{B}^p_{c,T}.$
In this formula, we may replace the variable $t$ with $t^T$ on both sides.
Furthermore, we may employ Proposition~\ref{prop3} for extending the
argument of the polar Mellin derivative from $(t,0)$ to an arbitary point
$(t, \alpha)\in \hs$. This leads us to
$$
\bigl(\widetilde{\Theta}_c f\bigr)(t,\alpha) \,=\,
\frac{4T}{\pi^2} \sum_{k=-\infty}^\infty \frac{(-1)^k}{(2k+1)^2}\, \,r_k^{c/T} f(tr_k^{1/T},\alpha)$$
for $f \in\mathscr{B}^p_{c,T}$ and all $(t,\alpha)\in\hs$,
which is identical with \zit{thm21}.

It remains to justify the statement on convergence. Since
\bgl{thm2p2a}
\frac{4}{\pi^2} \sum_{k=-\infty}^\infty \rez{(2k+1)^2}\,=\,1
\egl
and by property (iii) in the definition of $\mathscr{B}^p_{c,T}$,
there holds
$$ t^cr_k^{c/T} \abs{f\bigl(tr_k^{1/T}, \alpha\bigr)}\,\le\, 
C_f\, e^{T\abs{\alpha}}\quad (k\in\Z),$$
we readily see absolute and uniform convergence as asserted in the
theorem.
\beweisende


\section{Mellin--Bernstein inequality for polar Mellin derivatives}

As a consequence of Theorem~\ref{thm2}, we can establish a Bernstein
inequality for polar Mellin derivatives.

\begin{corollary}\label{cor1}
Let $f\in\mathscr{B}^p_{c,T}$, where $p\in [1, +\infty]$, $c\in\R$ and
$T>0$. Then
\bgl{cor11}
\li\|\bigl(\widetilde{\Theta}_cf\bigr)(\cdot, \theta)\ri\|_{X_c^p}\,\le\,
T\,\li\|f(\cdot, \theta)\ri\|_{X_c^p}
\egl
for any $\theta\in\R.$
\end{corollary}

\beweis
Let us write $r_k=e^{(k+1/2)\pi}$, where $k\in\Z$, for short.
Multiplying \zit{thm21} by $r^c$, taking moduli on both sides, applying
the triangular inequality on the right-hand side and using \zit{thm2p2},
we obtain
\begin{align}
r^c \abs{\bigl(\widetilde{\Theta}_cf\bigr)(r,\theta)} 
& \le
\frac{4T}{\pi^2}\sum_{k\in\Z} \rez{(2k+1)^2} \,r^cr_k^{c/T}
\abs{f(rr_k^{1/T}, \theta)} \label{cor1p1} \\[1ex]
&\le \frac{4T}{\pi^2} \sum_{k\in\Z} \rez{(2k+1)^2}\, \|f(\cdot,
\theta)\|_{X_c^\infty} \nonumber \\[1ex]
&= T\,\|f(\cdot, \theta)\|_{X_c^\infty} \nonumber
\end{align}
for any $(r, \theta)\in\hs$. This implies \zit{cor11} for $p=+\infty$.

Now let $1\le p< +\infty$ and let $R>0$. Observing that
$$  a^c \|f(a \,\cdot, \theta)\|_{X_c^p}\,=\, \|f( \cdot,\theta)\|_{X_c^p}$$
for $a>0$, we conclude from \zit{cor1p1} using the triangular inequality
for norms and B.\ Levi's theorem that
\begin{align*}
\li(\int_{1/R}^R r^{cp}\abs{\bigl(\widetilde{\Theta}_cf\bigr)(r,\theta)}^p
\frac{dr}{r}\ri)^{1/p}
&\le \li(\int_0^\infty r^{cp}\abs{\frac{4T}{\pi^2} \sum_{k\in\Z} 
\frac{r_k^{c/T}}{(2k+1)^2}
\,f(rr_k^{1/T},\theta)}^p \frac{dr}{r}\ri)^{1/p}\\[1ex]
&= \frac{4T}{\pi^2} \li\|\sum_{k\in\Z} \frac{r_k^{c/T}}{(2k+1)^2}\,
f(r_k^{1/T} \,\cdot, \theta)\ri\|_{X_c^p}\\[1ex]
& \le \frac{4T}{\pi^2} \sum_{k\in\Z} \frac{r_k^{c/T}}{(2k+1)^2}\,
\li\|f(r_k^{1/T} \,\cdot, \theta)\ri\|_{X_c^p}\\[1ex]
&= \frac{4T}{\pi^2} \sum_{k\in\Z} \rez{(2k+1)^2}\, \|f(\cdot,
\theta)\|_{X_c^p} \\[1ex]
&= T\,\|f(\cdot, \theta)\|_{X_c^p}\,.
\end{align*}
This guarantees that the $X_c^p$ norm  of
$\bigl(\widetilde{\Theta}_cf\bigr)(\cdot, \theta)$ exists. Letting
$R\to\infty$, we obtain \zit{cor11} for $p\in [1, \infty{[}$.
\beweisende

\begin{corollary}\label{cor2}
The Mellin-Bernstein space is invariant under polar Mellin
differentiation, that is, if $f\in\mathscr{B}^p_{c,T},$ where $p\in
[1,+\infty]$, $c\in\R$ and $T>0$, then $\widetilde{\Theta}_cf
\in\mathscr{B}^p_{c,T}$.
\end{corollary}

\beweis
We have to show that $\widetilde{\Theta}_c f$ satisfies properties
(i)--(iii) in the definition of Mellin--Bernstein spaces. Since for a
polar-analytic function $f$ there exist polar derivatives of arbitrary
order (see \cite[Theorem~5.3]{BBMS4}), it follows immediately that
$\widetilde{\Theta}_cf$ is polar-analytic on $\hs$, and so (i) is
satisfied. Property (ii) is guaranteed by Corollary~\ref{cor1}. Finally,
knowing that (iii) holds for $f$, we deduce from \zit{thm21} that
$$ r^c\abs{\bigl(\widetilde{\Theta}_cf\bigr)(r,\theta)}\,\le\,
\frac{4T}{\pi^2} \sum_{k\in\Z} \frac{C_f\,e^{T\abs{\theta}}}{(2k+1)^2}
\,=\, TC_f e^{T\abs{\theta}}.$$
Hence (iii) holds with $C_{\widetilde{\Theta}_cf} = T C_f$.
\beweisende

Combining Corollaries~\ref{cor1} and \ref{cor2}, we can apply \zit{cor11}
repeatedly and obtain
$$\li\|\bigl(\widetilde{\Theta}_c^kf\bigr)( \cdot , \theta)\ri\|_{X_c^p}
\,\le\, T^k\,\|f( \cdot, \theta)\|_{X_c^p}
\qquad (\theta\in\R, k\in\N).$$

\section{A Mellin analogue of Valiron's sampling formula}

In Fourier analysis, there exists a sampling formula of Valiron---sometimes also attributed to Tschakaloff
(see e.g. \cite{HIG}[p.~60], \cite{BSS})---that improves upon the classical sampling formula. It applies to a wider class of
functions and its series converges faster. The only price one has a pay for this improvement is that, apart from the
samples of the classical sampling formula, also the value of the derivative at zero is needed.

Here we want to establish an analogue of Valiron's formula in the Mellin setting. Our approach is very similar to the
derivation of an analogue of Boas' differentiation formula.

\begin{theorem}\label{thm6} Let $f\in \mathscr{B}_{c,T}^\infty$, where $c\in\R$ and $T\in\R^+$.
Then, for $r\in\R^+$,
\begin{align*}
 r^c f(r,0) = \sin(T\log r)&\biggl[\frac{(\widetilde{\Theta}_cf)(1,0)}{T} + \frac{f(1,0)}{T\log r}\\[1ex]
\quad&+ T\log r \sum_{k\in\Zz} \frac{(-1)^{k+1} e^{k\pi c/T} f(e^{k\pi/T},0)}{k\pi(k\pi -T\log r)}\biggr].
\end{align*}
The series converges absolutely and uniformly on compact subsets of $\R^+$.
\end{theorem}

\beweis
First we suppose that $T=1$. Let $t\in\R^+\setminus\{e^{k\pi}\,:\,k\in\Z\}$ be arbitrary but fixed and consider
the function
$$ F\::\: (r,\theta)\, \longmapsto\, \frac{f(r,\theta)}{(\log(r/t) +i\theta)(\log r +i\theta) \sin(\log r +i\theta)}\,.$$
It is polar-analytic on $\hs$ except for isolated singularities at the points where the denominator vanishes.
These points are $(t,0)$ and  $(e^{k\pi}, 0)$ for all $k\in\Z$. The point $(1,0)$, obtained for $k=0$, is a
logarithmic pole of order two while all other isolated singularities are logarithmic poles of order one.

Let $n\in\N$ be such that $e^{-n\pi}<t<e^{n\pi}$, and denote by $\partial\mathcal{R}_n$ the positively oriented rectangle
with vertices at the points $(e^{\pm(n+1/2)\pi}, \pm(n+1/2)\pi)$. By Theorem \ref{thm5}, we have
$$\rez{2\pi i} \int_{\partial\mathcal{R}_n} \bigl(re^{i\theta}\bigr)^{c-1} F(r,\theta) e^{i\theta}(dr +ird\theta)
=(\res_cF)(t,0) + \sum_{k=-n}^n (\res_cF)(e^{k\pi},0).$$

First we want to show that the integral on the left-hand side approaches zero as $n\to\infty$. Since $f$ was assumed to
belong to $\mathscr{B}_{c,1}^\infty$, we have
$$\abs{f(r,\theta)}\,\le\,C_f r^{-c} e^{\abs{\theta}}.$$
Furthermore it can be verified that on $\partial\mathcal{R}_n$
$$\abs{\sin(\log r+i\theta)}\,\ge\,\frac{e^{\abs{\theta}}}{3}\,.$$
Thus
\bgl{thm6p0}
\abs{r^{c-1}F(r,\theta)}\,\le\, \frac{3C_f}{r\bigl(\log^2(r/t) +\theta^2\bigr)^{1/2}\bigl(\log^2r+\theta^2\bigr)^{1/2}}\,.
\egl
For an  estimate of the contributions coming from the horizontal parts of $\partial\mathcal{R}_n$, we have to substitute $\theta=\pm(n+1/2)\pi$
in \zit{thm6p0}. Replacing in addition the squared
logarithms by zero, we obtain for the integrals under considerations the upper bound
$$\frac{3}{\pi} C_f \int_{e^{-(n+1/2)\pi}}^{e^{(n+1/2)\pi}} \rez{(n+1/2)^2 \pi^2} \frac{dr}{r},$$
which  approaches zero as $n\to\infty$.
For the contributions coming from the vertical parts we  have to substitute $r=e^{\pm(n+1/2)\pi}$ in \zit{thm6p0}. Replacing in addition $\theta^2$ by zero,
we obtain an upper bound for the corresponding integrals which again approaches zero as $n\to \infty.$
Hence the above residue formula yields
\bgl{thm6p1}
(\res_cF)(t,0)\,=\, -\,\sum_{k\in\Z} (\res_cF)(e^{k\pi},0).
\egl

Next we calculate the residues. It is easily seen that
$$(\res_cF)(t,0)\,=\, \frac{t^c f(t,0)}{\log t \sin(\log t)}.$$
For the residues on the right-hand side of \zit{thm6p1} we factor the sine in dependence of $k$ as
$$ \sin z\,=\, (z-k\pi) \psi_k(z) \qquad(z\in\C, \, k\in\Z).$$
Then
$$(\res_cF)(e^{k\pi},0) = \frac{e^{ck\pi} f(e^{k\pi},0)}{k\pi(k\pi-\log t) \psi_k(k\pi)}
=(-1)^k \frac{e^{ck\pi}f(e^{k\pi},0)}{k\pi(k\pi-\log t)} \quad (k\ne0).$$
Since in a neighborhood of $(1,0)$ we may write
$$F(r,\theta)\,=\, \rez{(\log r +i\theta)^2}\,\frac{f(r,\theta)}{(\log(r/t)+i\theta)\psi_0(\log r +i\theta)}
=: \frac{g(r,\theta)}{(\log r + i\theta)^2}\,,$$
we obtain
$$(\res_cF)(1,0)\,=\, (\widetilde{\Theta}_cg)(1,0)\,=\, (D_{\rm pol}g)(1,0)+ c g(1,0).$$
Noting that $\psi_0(z)=\sinc(z/\pi)$, we may express $g$ as
$$g(r,\theta)\,=\, \frac{f(r,\theta)}{\bigl(\log(re^{i\theta})-\log t\bigr)\sinc\bigl(\log(re^{i\theta})/\pi\bigr)}\,.$$
Now, via polar differentiation of the fraction on the right-hand side, we arrive at
$$(\res_cF)(1,0)\,=\,-\li[\frac{(\widetilde{\Theta}_cf)(1,0)}{\log t} + \frac{f(1,0)}{\log^2 t}\ri].$$

With these values for the residues, formula \zit{thm6p1} may be rewritten as
$$t^cf(t,0)= \sin(\log t)\li[(\widetilde{\Theta}cf)(1,0)+\frac{f(1,0)}{\log t}+
\log t \sum_{k\in\Zz}\frac{(-1)^{k+1}e^{ck\pi}f(e^{k\pi},0)}{k\pi(k\pi-\log t)}\ri].$$
It holds for $f\in\mathscr{B}_{c,1}^\infty$ and all $t\in\R^+$, even at the points $t=e^{k\pi}$
since
$$ \frac{\sin(\log t)}{k\pi -\log t}\,=\,(-1)^{k+1} \sinc\li(\frac{\log t}{\pi} -k\ri)$$
has a continuous continuation to the whole of $\R^+$.

With the help of Proposition \ref{prop2}, the achieved formula can be extended from $T=1$ to
abribrary $T\in\R^+.$ Then, substituting $t^{1/T}=r$, we obtain the formula stated in
Theorem~\ref{thm6}.

Noting that $e^{k\pi c/T}\abs{f(e^{k\pi/T},0)}\le C_f$ for all $k\in\Z$, we easily verify the assertion
on the convergence of the series.
\beweisende
\vskip0,3cm
\begin{Remark}\label{Rem1}
{\rm We could write the formula of Theorem~\ref{thm6} in terms of the $\lin$ function (see \cite{BJ2}), defined for $c \in \mathbb{R}$ and $x \in \mathbb{R}^+,$ by
$$\lin_c(x):= x^{-c} \mbox{sinc} (\log x),$$
with the continuous extension $\lin_c(1) := 1.$ 
Indeed, we have
\begin{align*}
f(r,0) &= \lin_{c\pi/T}(r^{T/\pi})\li[\log r \bigl(\widetilde{\Theta}_cf\bigr)(1,0) + f(1,0)\ri]\\[1ex]
& \quad + \log(r^{T/\pi}) \sum_{k\in\Zz}\frac{f(e^{k\pi/T},0)}{k} \lin_{c\pi/T}(e^{-k} r^{T/\pi}).
\end{align*}
This formula easily compares with the exponential sampling formula (\cite{BJ2}) and its right-hand side is defined for
all $r\in\R^+.$ On the other hand, it is very inconvenient when we want to estimate $f$ using property (iii)
of the Mellin-Bernstein space since all factors depending on $c$ are hidden in the $\lin$ function. An alternative
way would be to use $\lin_0$ only. This gives
 \begin{align*}
r^cf(r,0) &= \lin_0(r^{T/\pi})\li[\log r \bigl(\widetilde{\Theta}_cf\bigr)(1,0) + f(1,0)\ri]\\[1ex]
& \quad + \log(r^{T/\pi}) \sum_{k\in\Zz}\frac{e^{k\pi c/T}f(e^{k\pi/T},0)}{k} \lin_0(e^{-k} r^{T/\pi}).
\end{align*}
Now the factors depending on $c$ can go with $f$ as needed in (iii).}
\end{Remark}

Valiron's formula and its Mellin analogue will serve primarily for the reconstruction of a function $f$ from
samples. Among the latter, there is also a sample of the derivative of $f$. One may therefore think of interpreting this formula as a representation of the derivative of $f$ by samples of $f$.  In the following corollary,
we establish such a differentiation formula by skillful use of Theorem~\ref{thm6}. It has some advantages over the corresponding analogue of Boas' differentiation
formula. It seems that even a Fourier version of this result has not yet been mentioned in the literature.

\begin{corollary}\label{cor3}
Let $f\in\mathscr{B}_{c,T}^\infty$, where $c\in\R$ and $T>0$. Then
\begin{align*}
(\widetilde{\Theta}_cf)(r,\theta)\,=&\frac{T}{2}\li[e^{\pi c/(2T)}f(re^{\pi/(2T)},\theta)-e^{-\pi c/(2T)}f(re^{-\pi/(2T)},\theta)\ri]\\[1ex]
&\: +\frac{T}{\pi} \sum_{k\in\Zz}\frac{(-1)^k e^{k\pi c/T} f(re^{k\pi/T},\theta)}{k(4k^2-1)}
\end{align*}
for $(r,\theta)\in\hs$.
Multiplied by $r^c$, the series converges absolutely and uniformly on strips of bounded width parallel to the $r$-axis in $\hs$.
\end{corollary}

\beweis
By Proposition \ref{prop2}, the Mellin--Bernstein space is invariant under Mellin translations,  the formula of Theorem~\ref{thm6}
applies to
$$ g\::\: (r,\theta) \longmapsto t^c f(tr,\theta)\qquad (t>0).$$
Rewritten in terms of $f$, it takes the form
\begin{align*}
 r^c f(tr,0) = \sin(T\log r)&\biggl[\frac{(\widetilde{\Theta}_cf)(t,0)}{T} + \frac{f(t,0)}{T\log r}\\[1ex]
\quad&+ T\log r \sum_{k\in\Zz} \frac{(-1)^{k+1} e^{k\pi c/T} f(te^{k\pi/T},0)}{k\pi(k\pi -T\log r)}\biggr].
\end{align*}
Now we substitute $r=e^{\pm \pi/(2T)}$ and obtain the two equations:
\begin{align*}
e^{\pi c/(2T)} f(te^{\pi/(2T)},0)\,=&\frac{(\widetilde{\Theta}_cf)(t,0)}{T} + \frac{2}{\pi} f(t,0)\\[1ex]
&\: +\,\rez{2\pi} \sum_{k\in\Zz} \frac{(-1)^{k+1} e^{k\pi c/T} f(te^{k\pi/T},0)}{k(k-1/2)}
\end{align*}
and
\begin{align*}
e^{-\pi c/(2T)} f(te^{-\pi/(2T)},0)\,=&-\,\frac{(\widetilde{\Theta}cf)(t,0)}{T} + \frac{2}{\pi} f(t,0)\\[1ex]
&\: +\,\rez{2\pi} \sum_{k\in\Zz} \frac{(-1)^{k+1} e^{k\pi c/T} f(te^{k\pi/T},0)}{k(k+1/2)}
\end{align*}
Next, by subtraction of these equations, we find that
\begin{align*}
(\widetilde{\Theta}_cf)(t,0)\,=&\frac{T}{2}\li[e^{\pi c/(2T)} f(te^{\pi/(2T)},0)-e^{-\pi c/(2T)} f(te^{-\pi/(2T)},0)\ri]\\[1ex]
&\:+\,\frac{T}{\pi} \sum_{k\in\Zz} \frac{(-1)^k e^{k\pi c/T} f(t e^{k\pi/T},0)}{k(4k^2-1)}\,.
\end{align*}
Finally, employing Proposition \ref{prop3}, we can transform the argument of the polar Mellin derivative from $(t,0)$  to
an arbitrary point $(t,\alpha)\in\hs$. This leads us to
\begin{align*}
(\widetilde{\Theta}_cf)(t,\alpha)\,=&\frac{T}{2}\li[e^{\pi c/(2T)} f(te^{\pi/(2T)},\alpha)-e^{-\pi c/(2T)} f(te^{-\pi/(2T)},\alpha)\ri]\\[1ex]
&\:+\,\frac{T}{\pi} \sum_{k\in\Zz} \frac{(-1)^k e^{k\pi c/T} f(t e^{k\pi/T},\alpha)}{k(4k^2-1)}\,,
\end{align*}
which is the formula of the corollary with $(r,\theta)$ replaced by $(t,\alpha)$.
Recalling property (iii) in the definition of $\mathscr{B}_{c,T}^\infty$, we easily verify the statement on convergence.
\beweisende
\vskip0,3cm
\begin{Remark}\label{Rem2}
{\rm In order to write the differentiation formula of Corollary~\ref{cor3} in a simpler and more suggestive way, we introduce
the {\it central Mellin difference $\delta_{c,h}$ with increment $h>0$} by
$$ (\delta_{c,h}f)(r,\theta)\,:=\, h^c f(hr,\theta)- h^{-c}f(h^{-1}r,\theta).$$
\noindent
Then, setting $h:=e^{\pi/T}$, we may write the formula of Corollary~\ref{cor3} as
\bgl{new_diff}
(\widetilde{\Theta}_cf)(r,\theta)\,=\, T \li[\rez{2}\bigl(\delta_{c,h^{1/2}}f\bigr)(r,\theta)
+\rez{\pi} \sum_{k=1}^\infty \frac{(-1)^k}{k(4k^2-1)} \bigl(\delta_{c,h^k}f\bigr)(r,\theta)\ri].
\egl
Likewise the Mellin analogue of Boas' differentiation formula  may be expressed as
\bgl{Boas_diff}
(\widetilde{\Theta}_cf)(r,\theta)\,=\, \frac{4T}{\pi^2}
\sum_{k=0}^\infty \frac{(-1)^k}{(2k+1)^2} \bigl(\delta_{c,h^{k+1/2}}f\bigr)(r,\theta),
\egl
where again $h=e^{\pi/T}$.
We observe that the summation in \zit{new_diff} has coefficients decaying like $\mathcal{O}(k^{-3})$ as $k\to \infty$ while those of
\zit{Boas_diff} decay like $\mathcal{O}(k^{-2})$ only.}
\end{Remark}
\smallskip
\begin{Remark}\label{Rem3}
{\rm
For sake of completeness, let  us have a quick look at the classical Fourier  case. For $f$ belonging to the (Fourier) Bernstein space
$B_w^\infty$, Valiron's sampling formula states that
$$f(t)=\sin(wt)\li[\frac{f'(0)}{w} + \frac{f(0)}{wt} + wt \sum_{k\in\Zz}\frac{(-1)^k f(k\pi/w)}{k\pi(wt-k\pi)}\ri]$$
for $t\in\R$. Now we may imitate the proof of Corollary~\ref{cor3} in a Fourier version, that is, we apply this formula
to $f(\cdot + x)$, substitute in the resulting formula $t=\pm \pi/(2w)$ and subtract the two resulting equations.
Thus we arrive at the (new) differentiation formula
\begin{align*}
f'(x)&= \frac{w}{2}\li[f(x+\pi/(2w))-f(x-\pi/(2w))\ri] + \frac{w}{\pi}\sum_{k\in\Zz}\frac{(-1)^k f(x+k\pi/w)}{k(4k^2-1)}\\[1ex]
&= \frac{w}{2}\li[f(x+\pi/(2w))-f(x-\pi/(2w))\ri]\\[1ex]
& \quad + \frac{w}{\pi}\sum_{k=1}^\infty\frac{(-1)^k}{k(4k^2-1)}
\li[f(x+k\pi/w)-f(x-k\pi/w)\ri].
\end{align*}
Thus the derivative is represented by an alternating series of central differences.
Of course, in this formula we can replace the real variable $x$ also by the complex variable $z$.
}
\end{Remark}
\vskip0,3cm
We conclude the section with an interesting consequence of Theorem \ref{thm6}, which states an identity theorem in Mellin--Bernstein spaces.

\begin{corollary}\label{cor1i}
Let $f\in \mathscr{B}_{c,T}^\infty$, where $c\in\R$ and $T>0$. If $f$ vanishes on the set $\bigl\{(e^{k\pi/T},0)\,:\, k\in\Z\bigr\}$,
then
$$ f(r,\theta)\,\equiv\, a\bigl(re^{i\theta}\bigr)^{-c} \sin\bigl(T(\log r +i\theta)\bigr)$$
for some $a\in\C.$
\end{corollary}

\beweis
It follows from the Mellin analogue of Valiron's sampling formula (Theorem \ref{thm6}) that
$$ f(r,0)\,=\, a r^{-c} \sin(T\log r)$$
for all $r>0$, where $a=(\wT_cf)(1,0)/T.$

Now consider the function $g\,:\,\hs\to\C$ defined by
$$g(r,\theta)\,=\, a \bigl(re^{i\theta}\bigr)^{-c} \sin\bigl(T(\log r+i\theta)\bigr).$$
Obviously it belongs to $\mathscr{B}_{c,T}^\infty$ and $g(r,0)=f(r,0)$ for all $r>0$. Thus $\varphi:=f-g$ also belongs to
$\mathscr{B}_{c,T}^\infty$ and $\varphi(r,0)=0$ for all $r>0.$
By Theorem~\ref{thm1i}, it follows that $\varphi(r,\theta)\equiv 0$ and so $f(r,\theta)\equiv g(r,\theta)$ on $\hs$.
\beweisende

As a consequence of Corollary~\ref{cor1i}, we obtain the following result:
\begin{corollary}\label{cor2i}
Let $f\in \mathscr{B}_{c,T}^p$, where $c\in\R$, $T>0$ and $p\in [1, \infty{[}$. If $f$ vanishes on the set $\bigl\{(e^{k\pi/T},0)\,:\, k\in\Z\bigr\}$,
then $f$ is identically zero.
\end{corollary}

\beweis
Since $\mathscr{B}_{c,T}^p\subset \mathscr{B}_{c,T}^\infty$, Corollary~\ref{cor1i} applies to $f$ and yields that
$$ f(r,\theta)\,\equiv\, a\bigl(re^{i\theta}\bigr)^{-c} \sin\bigl(T(\log r +i\theta)\bigr)$$
for some $a\in\C.$ The membership of $f$ in $\mathscr{B}_{c,T}^p$ requires that $f(\cdot, 0)\in X_c^p$. This is
possible if and only if $a=0$.
\beweisende

\section{Computational Aspects}

For practical use of  formulae \zit{new_diff} and \zit{Boas_diff} we have to truncate the series.
As an approximation of the polar Mellin derivative by $2n$ samples, we have
\bgl{new_diff_n}
(\widetilde{\Theta}_cf)(r,\theta)\,=\, T \li[\rez{2}\bigl(\delta_{c,h^{1/2}}f\bigr)(r,\theta)
+\rez{\pi} \sum_{k=1}^{n-1} \frac{(-1)^k}{k(4k^2-1)} \bigl(\delta_{c,h^k}f\bigr)(r,\theta)\ri]
+ E_{1,n}(f)
\egl
and
\bgl{Boas_diff_n}
(\widetilde{\Theta}_cf)(r,\theta)\,=\, \frac{4T}{\pi^2}
\sum_{k=0}^{n-1} \frac{(-1)^k}{(2k+1)^2} \bigl(\delta_{c,h^{k+1/2}}f\bigr)(r,\theta)
+ E_{2,n}(f)
\egl
with error terms $E_{1,n}(f)$ and $E_{2,n}(f)$, respectively.
Employing property (iii) of the Mellin--Bernstein space, we see that
$$\abs{ (\delta_{c,t}f)(r,\theta)}\,\le\, 2C_f r^{-c} e^{T\abs{\theta}}$$
for any increment $t>0$.
Therefore we can estimate the errors as
$$ \abs{E_{1,n}(f)}\,\le\, \frac{2C_f}{\pi}\,r^{-c} Te^{T\abs{\theta}} \sum_{k=n}^\infty \rez{k(4k^2-1)}$$
and
$$ \abs{E_{2,n}(f)}\,\le\, \frac{8C_f}{\pi^2}\,r^{-c} Te^{T\abs{\theta}} \sum_{k=n}^\infty \rez{(2k+1)^2}.$$
For estimating the series on the right-hand side, we may use integral comparison and obtain:
\begin{align*}
\sum_{k=n}^\infty \rez{k(4k^2-1)} &\le \int_{n-1}^\infty \frac{dx}{x(4x^2-1)} \\[1ex]
&= -\,\rez{2} \log\li(1- \rez{4(n-1)^2}\ri)\\[1ex]
&= \rez{2} \sum_{j=1}^\infty \rez{j[4(n-1)^2]^j}\\[1ex]
&\le \rez{2} \sum_{j=1}^\infty \rez{[4(n-1)^2]^j}\,
= \rez{8(n-1)^2-2}
\end{align*}
and
$$ \sum_{k=n}^\infty \rez{(2k+1)^2}\,\le\, \int_{n-1}^\infty \frac{dx}{(2x+1)^2}\,=\,\rez{4n-2}\,.$$
Thus we arrive at the error bounds
$$ \abs{E_{1,n}(f)}\,\le\, \frac{C_f r^{-c} T e^{T\abs{\theta}}}{\pi[4(n-1)^2-1]}$$
and
$$ \abs{E_{2,n}(f)}\,\le\, \frac{4C_f r^{-c} T e^{T\abs{\theta}}}{\pi^2(2n-1)}\,.$$

\section{ A short biography of Domenico Candeloro 1951--2019}

Domenico Candeloro was born October 18, 1951 in Udine, northern Italy, while his parents descended from Abruzzo, a large region in central Italy, about 90 km east from Rome. He obtained the laurea degree in Mathematics with highest distinction under Professor Calogero Vinti in July 1974, his thesis being in the field of Calculus of Variations. His interests soon expanded to Real Analysis and Measure Theory, and the fields above mentioned.

Domenico belonged to Calogero Vinti second set of students, together with Patrizia Pucci, Anna Salvadori, all three being 67 years old. The first set consisted of Prof. Mauro Boni (died 2010 at 67), Primo Brandi, Marcello Ragni and Candida Gori Cocchieri all between 71 to 73 years old. The third set consisted of Carlo Bardaro, Anna Martellotti and Rita Ceppitelli, all 64 years old. Students of Calogero’s students, all under 60, include Gianluca Vinti, Anna Rita Sambucini, Ilaria Mantellini, Antonio Boccuto, Luca Zampogni, Laura Angeloni and Danilo Costarelli.

Through Professor Calogero Vinti, Domenico himself was an academic great-grandson of Leonida Tonelli, in turn a great-great grandson of Cesare Arzel\'a, whose teacher was Enrico Betti, a great lineage indeed.

Since Mimmo became Full Professor in 1997, but also as associate Professor, he taught Probability Theory, Functional Analysis and Stochastic Processes at Perugia, and was thesis advisor of many students. One often regards him as the most brilliant student of the Vinti School.

One of his more recent papers, entitled {\it Vitali-type theorems for filter convergence related to vector lattice-valued modular and applications to stochastic processes}, published in Journal of Mathematical Analysis and Applications, in 2014, and written in collaboration with Anna Rita Sambucini and Antonio Boccuto, was awarded with the prestigious 2014 JMAA Ames Award. The results contained in this paper establish deep connections with apparently different areas of mathematics, giving a link between Brownian motions, approximation theory and Mellin-type convolution operators.

He collaborated with numearous colleagues outside the University of Perugia, thus with J.K. Brooks (Florida, USA), A. Volcic (Trieste), K.P.S. Bhaskara Rao (India), G. Letta (Pisa), K. Musial and  E. Kubinska (Poland), B. Riecan and R. Mesiar (Slovacchia), A.C. Gravilut and A. Croitoru (Romania), C.A. Labuschagne (South Africa), L. Di Piazza and V. Marraffa (Palermo), the many collaborations internally included P. Pucci, Anna Martellotti, Anna Rita Sambucini, Antonio Boccuto, Carlo Bardaro. He was the author of more than 75 papers in established journals. 
He was invited by many universities to hold seminars and conferences, among them, Florida State University at Gainesville, Universities in Oxford, Reading, London (England), University of Bratislava, Nitra, Banska, Bystrica (Slovacchia). He was among the organizers of several scientific initiatives and the coordinator of many research projects.
\vskip0,3cm
\noindent
A personal remembrance of Paul Butzer:\\
 {\it Mimmo, Mimmo, please help me, I am stuck. Mimmo, Mimmo, I think that I have found a new theorem, another voice noted; Mimmo looked at the proof, with the words} I feel it is an easy application of a classical result. {\it It is these calls for assistance which drew my attention to the mathematician himself, Domenico Candeloro, “Mimmo” being the familiar name his close friends call him. It is his expertise in many areas, Calculus of Variations, Measure Theory, Real Analysis, Functional Analysis, Probability Theory, which place him in a central role in Perugia, one who is always willing to help.

At a few meetings with Mimmo, during my first years at Perugia from 1990 onward and my participations in the conference of the Research Group} Real Analysis and Measure Theory, {\it conducted by Calogero Vinti and Colleagues in Capri, Ischia, Maiori, and Grado, together with my friends Carlo, Gianluca and Anna Rita, Mimmo was in the company of his Spouse, Prof. Doretta Vivona. He seemed happy that she carried on the conversation (in English). Both were very friendly, Mimmo having an engaging and modest personality}.

\section*{Acknowledgments}
{\small  The authors express their gratitude to Dr. Cecilia Bracuto for her help in drafting some graphic representations.  Carlo Bardaro and Ilaria Mantellini have been partially supported by the ``Gruppo Nazionale per l'Analisi Matematica e Applicazioni (GNAMPA)'' of the
``Instituto di Alta Matematica (INDAM)'' as well as by the projects ``Ricerca di Base 2017 of University of Perugia (title: Misura, Integrazione,
Approssimazione e loro Applicazioni)'' and ``Progetto Fondazione Cassa di Risparmio cod. nr. 2018.0419.021 (title: Metodi e Processi di Intelligenza
artificiale per lo sviluppo di una banca di immagini mediche per fini diagnostici (B.I.M.))''.}


\end{document}